\documentclass[a4paper]{article}
\usepackage{amsthm,amsmath,amssymb}
\usepackage{color}
\newtheorem{teor}{Theorem}[section]
\newtheorem{defin}[teor]{Definition}
\newtheorem{lemm}[teor]{Lemma}
\newtheorem{osse}[teor]{Remark}
\newtheorem{prop}[teor]{Proposition}
\newtheorem{defi}[teor]{Definition}
\newtheorem{coro}[teor]{Corollary}
\newtheorem{prob}[teor]{Problem}
\newtheorem{assu}[teor]{Assumption}

\newcommand{\bele}{\begin{lemm}\begin{sl}}
\newcommand{\enle}{\end{sl}\end{lemm}}
\newcommand{\bedef}{\begin{defi}\begin{sl}}
\newcommand{\eddef}{\end{sl}\end{defi}}
\newcommand{\bete}{\begin{teor}\begin{sl}}
\newcommand{\ente}{\end{sl}\end{teor}}
\newcommand{\beos}{\begin{osse}\begin{rm}}
\newcommand{\eddos}{\end{rm}\end{osse}}
\newcommand{\beas}{\begin{assu}\begin{rm}}
\newcommand{\eddas}{\end{rm}\end{assu}}
\newcommand{\bepr}{\begin{prop}\begin{sl}}
\newcommand{\empr}{\end{sl}\end{prop}}
\newcommand{\bepro}{\begin{prob}\begin{rm}}
\newcommand{\empro}{\end{rm}\end{prob}}
\newcommand{\bede}{\begin{defin}\begin{sl}}
\newcommand{\edde}{\end{sl}\end{defin}}
\newcommand{\beco}{\begin{coro}\begin{sl}}
\newcommand{\enco}{\end{sl}\end{coro}}

\newcommand{\weak}{strong }

\newcommand{\quext}{\quad\text}
\newcommand{\qquext}{\qquad\text}

\newcommand{\RR}{\mathbb{R}}

\newcommand{\CC}{\mathbb{C}}

\DeclareMathOperator*{\esssup}{ess\,sup}

\newcommand{\beeq}[1]{\begin{equation}\label{#1}}
\newcommand{\eddeq}{\end{equation}}

\newcommand{\beeqa}[1]{\begin{eqnarray}\label{#1}}
\newcommand{\eddeqa}{\end{eqnarray}}

\newcommand{\beal}[1]{\begin{align}\label{#1}}
\newcommand{\eddal}{\end{align}}

\newcommand{\bespl}[1]{\begin{split}\label{#1}}
\newcommand{\edspl}{\end{split}}

\newcommand{\bega}[1]{\begin{gather}\label{#1}}
\newcommand{\edga}{\end{gather}}

\newcommand{\beeqax}{\begin{eqnarray*}}
\newcommand{\eddeqax}{\end{eqnarray*}}

\def\qed{\ifmmode 
  \else \leavevmode\unskip\penalty9999 \hbox{}\nobreak\hfill
  \fi
  \quad\hbox{\hskip.5em\vrule width.4em height.6em depth.05em\hskip.1em}}
\def\endproofsym{\qed}

\def\endnobox{\def\endproofsym{}\end{proof}\def\endproofsym{\qed}}

\newcommand{\no}{\nonumber}

\newcommand{\beeqao}{\begin{eqnarray}\no}
\newcommand{\bealo}{\begin{align}\no}
\newcommand{\besplo}{\begin{split}\no}
\newcommand{\begao}{\begin{gather}\no}

\newcommand{\duav}[1]{\langle{#1}\rangle}

\newcommand{\+}{\hspace{1pt}}

\newcommand{\itt}{\int_0^t}

\newcommand{\io}{\int_\Omega}

\newcommand{\iTT}{\int_0^T}

\newcommand{\epsi}{\varepsilon}

\newcommand{\lla}{_{\lambda}}

\newcommand{\dd}{_\delta}
\newcommand{\ddm}{_{\delta,m}}

\newcommand{\OO}{_{\Omega}}

\def\R{\mathbb R}

\newcommand{\bn}{\boldsymbol{n}}

\newcommand{\dn}{\partial_{\bn}}
\newcommand{\fhi}{\varphi}

\newcommand{\vu}{\boldsymbol{u}}

\newcommand{\bbg}{\boldsymbol{g}}

\newcommand{\lhs}{left-hand side}
\newcommand{\rhs}{right-hand side}

\DeclareMathOperator{\dive}{div}
\DeclareMathOperator{\deriv}{d}

\DeclareMathOperator{\exte}{ext}

\newcommand{\Bext}{B_{\+\exte}}

\let\TeXchi\chi
\def\chi{{\setbox0 \hbox{\mathsurround0pt
$\TeXchi$}\hbox{\raise\dp0 \copy0 }}}

\newcommand{\calX}{{\mathcal X}}

\newcommand{\calA}{{\mathcal A}}

\newcommand{\calE}{{\mathcal E}}

\newcommand{\calV}{{\mathcal V}}

\newcommand{\barO}{\overline{\Omega}}

\newcommand{\baru}{\overline{u}}

\newcommand{\dit}{\deriv\!t}

\newcommand{\dir}{\deriv\!r}

\newcommand{\ddt}{\frac{\deriv\!{}}{\dit}}

\newcommand{\Lpsot}{L^\psi(\Omega)}

\newcommand{\Lpsq}{L^\psi(Q)}

\newcommand{\deo}{\partial_\Omega}

\newcommand{\deQ}{\partial_Q}

\newenvironment{giuliorev}{\color{red}}{\color{black}}
\newcommand{\III}{\begin{giuliorev}}
\newcommand{\EEE}{\end{giuliorev}}

\numberwithin{equation}{section}
\begin{document}

\title{Local well-posedness for Fr\'emond's model\\of complete damage in elastic solids}

\author{Goro Akagi\\
Mathematical Institute and Graduate School of Science, Tohoku University,\\
6-3 Aoba, Aramaki, Aoba-ku, Sendai 980-8578 Japan\\
E-mail: {\tt goro.akagi@tohoku.ac.jp}\\
\and
Giulio Schimperna\\
Dipartimento di Matematica, Universit\`a di Pavia,\\
Via Ferrata~5, I-27100 Pavia, Italy\\
E-mail: {\tt giusch04@unipv.it}
}


\maketitle
\begin{abstract}
 We consider a model for the evolution of damage in elastic materials originally 
 proposed by Michel Fr\'emond. For the corresponding PDE system we prove existence and
 uniqueness of a local in time \weak solution. The main novelty of our result 
 stands in the fact that, differently from previous contributions,
 we assume no occurrence of any type of regularizing terms.
\end{abstract}


\vspace{2mm}

\noindent {\bf AMS (MOS) subject clas\-si\-fi\-ca\-tion:}%
~~35K55, 35K86, 35M33, 35Q74, 74R05.

\vspace{2mm}


\section{Introduction}
\label{sec:intro}

We consider a basic model for the evolution of damage in an elastic material
subject to an external load under the approach originally proposed by Fr\'emond and
coauthors in a number of papers \cite{FKS,FKNS,FN} (see also the monographs
\cite{Fre1,Fre2,Lem} for a general presentation of related models as well
as a detailed mechanical background).

We will give here an overview of the model in its generality; we notice
however from the very beginning that, in order to reduce technical complications, 
a simplified formulation will be addressed for the purpose of a mathematical
analysis. Let us consider a smooth and bounded domain
$\Omega\subset \RR^3$ occupied by the elastic medium over some given reference 
time interval $(0,T)$. The material is subject
to an external load $\bbg$ leading to elastic deformations represented
by means of the displacement variable $\vu$. As a response to deformations, the 
material undertakes elastic stresses that are a source of {\sl damage}.
At a microscopic level, this phenomenon can be thought as a progressive failure
of elastic bonds; as a consequence, the material loses stiffness and 
micro-cracks tend to develop.

A description of the progression of damage at the microscopic level is however
very difficult, especially because the micro-breaks are very small compared
to the scale of macroscopic displacements. For this reason, in this type of 
{\sl continuum models}\/ the damage is rather described by means of a {\sl macroscopic}\/ 
variable $z$, i.e., an order parameter that represents the locally 
averaged evolution of damage at any point $x\in \Omega$ and 
$t\in (0,T)$. For simplicity $z$ is normalized in such a way that, 
for $z=1$, the material is completely integer, i.e., no damage has yet occurred,
whereas for $z=0$ all the elastic bonds have been broken. We speak
then of {\sl complete damage}\/ at that point, meaning that the material
has completely lost its elastic properties and a (macroscopic) fracture 
has occurred. According to such an interpretation,
the values of $z$ below $z=0$, as well as those above $z=1$
have no physical significance and should be somehow penalized in the mathematical
formulation of the model. 

We will assume a quasi-static regime; namely, the damage process occurs at
a much slower scale compared to the elastic response, which can thus be represented
by an {\sl elliptic}\/ equation of the form
\begin{equation}\label{eq:u:compl}
  - ( \CC_{ijkl}(z) \epsi(\vu)_{kl} )_{,j} = g_i, \qquext{in }\,(0,T)\times \Omega.
\end{equation}
Here $\epsi(\vu) = ( \nabla \vu + (\nabla \vu)^t )/2$ is the  strain tensor, 
$\bbg=(g_i)$ represents the action of the (given) external forces,
and the elastic tensor $\CC$ may be assumed to satisfy proper symmetry and ellipticity
conditions and to degenerate as $z=0$ (the precise hypotheses will be presented below). 
Here and below we are assuming Einstein's convention for summation over repeated indices.
It is worth noting that {\sl dynamical}\/ models for damage evolution are also
significant and have been studied mathematically in a number of contributions.
We may quote, with no claim of completeness, \cite{BSS,FKS,GO,HK2,HKRR} 
(see also the references therein) for models including inertial and/or viscosity effects.

Relation \eqref{eq:u:compl} is complemented with the following parabolic equation describing the evolution of the damage variable $z$:
\begin{equation}\label{eq:z:compl}
   \alpha(z_t) + \delta_1 z_t - \delta_2 \Delta z + f'(z) 
    \ni w - \frac12 \CC'_{ijkl}(z) \epsi(\vu)_{kl} \epsi(\vu)_{ij}, \qquext{in }\,(0,T)\times \Omega,
\end{equation}
Here, $\alpha = \partial I_{(-\infty,0]}$, i.e., the subdifferential of the 
{\sl indicator function}\/ of the interval $(-\infty,0]$. We refer the reader to
the monographs \cite{Ba,Br} for the underlying background material from convex analysis.
Here we just recall that $\alpha$ is a multivalued mapping; indeed, we have
$\alpha(0) = [0, + \infty)$, $\alpha(z) = \{0\}$, for $z<0$ and $\alpha(z) = \emptyset$ 
for $z>0$. This motivates the occurrence of the inclusion sign in~\eqref{eq:z:compl}.
The presence of $\alpha$ is aimed at enforcing the 
{\sl irreversibility}\/ (or {\sl unidirectionality}) constraint on the evolution
of $z$. Namely, any solution must satisfy $z_t\le 0$, which means that
once some amount of damage has been created, it cannot be repaired. 
Note that this fact implies in turn that, once $z_0\le 1$, then $z$ can never 
exceed $1$ at any point in the evolution, implying that the unphysical
states $z>1$ are automatically excluded.
Irreversibility is a reasonable physical ansatz in many real world applications; 
on the other hand it is worth observing that also {\sl reversible}\/ models 
(i.e., such that the broken bonds may be at least partially restored) are significant 
and have been extensively studied in the literature (see, e.g., \cite{BFS} and the
references quoted there). It is also worth noticing that \eqref{eq:z:compl}
subsumes a {\sl rate-dependent}\/ evolution of $z$; {\sl rate-independent}\/
damage models are equally interesting and have been addressed in several
works (see, e.g., \cite{FKStef,KRZ,Mi,MR} and the references therein).

The coefficients $\delta_1,\delta_2>0$ in \eqref{eq:z:compl} are 
related to the time scale of the damaging process (the smaller $\delta_1$ the faster
it occurs) and to the ``thickness'' of the (diffuse) interface between damaged and 
sound areas (which depends on the scale length of the micro-breaks
and goes like $\delta_2^{1/2}$). The positive constant $w >0$ on 
the \rhs\ has the significance of a threshold:  let us explain this fact by
assuming $f\equiv 0$, which, physically speaking, can be seen 
as the ``model case''. In this situation, if the forcing term 
$\CC'_{ijkl}(z) \epsi(\vu)_{kl} \epsi(\vu)_{ij}$  does not exceed $2w$, 
the \rhs\ of \eqref{eq:z:compl} is positive, which basically indicates that no damage is being created.
In the converse situation, i.e.~in presence of large deformation gradients,
a source of damage occurs. In the case $f\not\equiv 0$, this 
damaging effect can be thought to vary a little depending
on the actual value of $z$; nevertheless one expects that,
in practice, $f'(z)$ is small compared to $w$. Hence, if we set
$\psi'(r) = f'(r) - w$ (as we will do in the sequel), we expect
in particular $\psi'$ be strictly negative or, in other words,
the {\sl configuration potential}\/ $\psi$ to be concave,
meaning that, in some measure, the body tends to oppose 
resistance to the damaging effects which, as said, will occur only 
if the elastic stresses are large. 

In order to present our mathematical results, let us assume for simplicity
$\bbg$ independent of time and take homogeneous Dirichlet boundary
conditions for $\vu$ and no-flux (i.e., homogeneous Neumann)
boundary conditions for $z$. Moreover, let us assume (at least)
the symmetry property $\CC_{ijkl} = \CC_{klij}$.
Then, testing \eqref{eq:u:compl} by $\vu_t$ and \eqref{eq:z:compl} 
by $z_t$ and integrating over $\Omega$ permits us to (formally) 
deduce the {\sl energy equality}
\begin{equation}\label{eq:E:compl}
  \ddt \calE(t) + \delta_1 \| z_t \|_{L^2(\Omega)}^2 = 0
\end{equation}
with the energy functional
\begin{equation}\label{E:compl}
   \calE(t) = \io \Big( \frac12 \CC_{ijkl}(z) \epsi(\vu)_{kl} \epsi(\vu)_{ij} 
    - \bbg \cdot \vu
    + \frac{\delta_2}2 |\nabla z |^2 + f(z) - w z \Big),
\end{equation}
where it is worth noting that the product between $z_t$ and $\alpha(z_t)$ 
is a.e.~equal to $0$, in view of the fact that $\alpha(z_t)$ (or, to be precise, 
any element of such a set) may be different from $0$ only when $z_t=0$.
The energy relation \eqref{eq:E:compl} 
is the basic source of the a priori estimates needed for
attempting a mathematical analysis of system \eqref{eq:u:compl}-\eqref{eq:z:compl}.

On the other hand, there are several reasons why the information provided by
the above relation is not sufficient in order to obtain a satisfactory mathematical
result. An important point stands of course in the fact that, even  if the body 
is completely integer at the beginning (i.e., $z_0\equiv 1$ in $\Omega$), 
it is expected that after some time, due to progression of damage, $z$ becomes $0$ 
at some point $x\in \Omega$. In such a situation, the elastic tensor $\CC(z)$ degenerates
and the energy $\calE$ is no longer coercive. Consequently, it becomes 
impossible to control the quadratic term in $\epsi(\vu)$ on 
the \rhs\ of \eqref{eq:z:compl} and the model somehow loses significance.
This is an intrinsic feature of this system (and of related ones) and,
actually, for such models of {\sl complete} damage, 
it seems natural to look for {\sl local in time solutions}, namely those defined on
a ``small'' time interval $(0,T_0)$ with possibly $T_0<T$, where degeneration does not occur.
This type of {\sl local existence}\/ result is what is proved in 
several related papers (see, e.g., \cite{BS,FKS,FKNS}) and will also be 
the object of the present note. Indeed, it seems that the description of 
{\sl complete}\/ damaging of the material, i.e., of what happens 
after the onset of some macroscopic fracture,
 requires a different modeling approach, see, e.g., \cite{BMR,HK,Mi}.

There is, however, a second relevant difficulty; indeed, in order to 
prevent degeneration of $z$ at least in a short time interval $(0,T_0)$,
one needs a quantitative estimate of the form
\begin{equation}\label{add:est}
  \| z \|_{L^\infty(0,T_0;X(\Omega))} \le c,
\end{equation}
where $T_0>0$ may depend on the prescribed data and $X=X(\Omega)$ is a Banach space such that $X\subset C^0(\barO)$ with 
continuous embedding. This corresponds to a (local) control of $z$ in the {\sl uniform}\/
norm, in such a way that degeneration cannot occur at any point in the short time span.
On the other hand, if the energy has the expression \eqref{E:compl}, an 
estimate like \eqref{add:est} follows directly from \eqref{eq:E:compl} only
in space dimension one (this is, indeed, the spirit of the pioneering results
proved in \cite{FKS,FKNS}), whereas, in the present three-dimensional
setting, \eqref{add:est} may be obtained only by performing 
higher regularity estimates. Here, however, two additional difficulties
arise: (i) the combined occurrence in \eqref{eq:z:compl} 
of the  {\sl nonsmooth}\/ function $\alpha$ and of the  {\sl quadratic gradient term}\/
on the \rhs, and (ii) the poor regularity of $\vu$ provided by the 
{\sl elliptic}\/ equation \eqref{eq:u:compl} characterized by a $z$-dependent
(hence nonsmooth) diffusion coefficient. For these reasons, at least up 
to our knowledge, local existence has been obtained so far only in presence of 
additional smoothing terms. Actually, common regularizations considered in the literature
are: {\sl viscoelastic}\/ (rather than purely elastic)
behavior for $\vu$ \cite{BBR,BSS,GO,HKRR}, 
presence of inertial effects in \eqref{eq:u:compl}
\cite{BSS,GO,HK2,HKRR}, and replacement of the Laplacian in \eqref{eq:z:compl}
by a more regularizing operator like the fractional Laplacian $(-\Delta)^s$
with suitable $s>1$ \cite{KRZ} or the $p$-Laplacian $-\Delta_p$ with
suitable $p>2$ \cite{HK,HK2,HKRR}.

In this work, we will consider the ``original'' system \eqref{eq:u:compl}-\eqref{eq:z:compl} 
with {\sl no occurrence of any regularizing term}. We will actually prove
that an estimate of the form \eqref{add:est} can be obtained also in such a 
setting, so filling the gap of a long-standing regularity problem. Our argument
is based on a more careful control of the $L^\infty$-, $H^1$- and $H^2$- norms of the difference
between $z(t)$ at $t>0$ and the initial datum $z_0$ in terms of the parameters of the system. 
As an outcome of our procedure, we will be able to prove existence and uniqueness 
of \weak solutions to the initial-boundary value problem for system 
\eqref{eq:u:compl}-\eqref{eq:z:compl} on a time span $(0,T_0)$, with 
$T_0$ explicitly computable in terms of the data, where $z$ 
does not degenerate to $0$ at any point. 

In order to avoid unessential technicalities, proceeding in the spirit of \cite{BS,BSS}
we will actually consider a simplified version of the model, where the displacement $\vu$ is replaced
by a scalar variable $u$ and some quantities and parameters are normalized.
We point out that these simplifications are not restrictive and are taken
only for the sake of clarity. Indeed, our results could be easily extended to the 
``original'' system \eqref{eq:u:compl}-\eqref{eq:z:compl} by applying some more 
or less standard tools (like, e.g., Korn's inequality) and doing a little more technical
work.

\medskip

The paper is organized as follows. In the next section, we provide a 
detailed presentation of our assumptions and state our main result.
The a-priori estimates that are at the core of the proof are given in
the subsequent Section~\ref{sec:proofs}. Finally, a possible regularization
of the system compatible with the a-priori estimates  
is sketched in the final Section~\ref{sec:appro}, where a number of 
additional comments are also given.


\section{Main result}
\label{sec:main}

First of all, we introduce a simplified version 
of system \eqref{eq:u:compl}-\eqref{eq:z:compl}. As said,
we replace the vector-valued displacement $\vu$ by a scalar one $u$, and
correspondingly assume that the elasticity tensor $\CC(z)$ is replaced by
a scalar function $c(z)$. Moreover, in order to take the simplest example
of a strictly positive function that degenerates at $0$ we just choose $c(z)=z$.
We also normalize the parameters $\delta_1$, $\delta_2$ to~$1$ and 
incorporate the positive constant $w$ into the function $f'$
so introducing a new configuration potential $\psi(r) = f(r) - wr$.
With these choices, system \eqref{eq:u:compl}-\eqref{eq:z:compl} reduces to
\begin{alignat}{4}\label{eq:u:strong}
 &-\dive(z\nabla u)=g, &&\qquext{in }\,(0,T)\times \Omega,\\
 \label{eq:z:strong}
 & \alpha(z_t) + z_t - \Delta z + \psi'(z) \ni - \frac12 | \nabla u |^2,
    &&\qquext{in }\,(0,T)\times \Omega.
\end{alignat}
The above equations are complemented with the boundary conditions
(which are a rather standard choice for this class of models) 
\begin{equation}\label{no:flux}
  u = \dn z = 0, \qquext{in }\,(0,T)\times \Gamma,
\end{equation}
where $\Gamma=\partial \Omega$,  $\dn = \bn \cdot \nabla$ and $\bn$ 
denotes the outer unit normal
vector to $\Gamma$. System \eqref{eq:u:strong}-\eqref{eq:z:strong}
is stated over an assigned reference interval $(0,T)$; however, 
as said, we will prove existence on a possibly smaller interval $(0,T_0)$. 
Finally, we assume the initial condition 
\begin{equation}\label{init}
  z|_{t=0} = z_0, \qquext{in }\Omega.
\end{equation}

\smallskip

In order to fix a concept of \weak solution and formulate our related existence 
result, we need to introduce some preparatory material.
Letting $\Omega$ be a smooth bounded domain of $\RR^3$, we 
set $H := L^2(\Omega)$, $V := H^1(\Omega)$ and $V_0 := H^1_0(\Omega)$. 
We will often write $H$ in place of $H \times H \times H$ 
(with similar notation for other spaces),
in case vector-valued functions are considered. We denote
by $(\cdot,\cdot)$ the standard scalar product of~$H$ and by 
$\| \cdot \|$ the associated Hilbert norm.  Moreover, we equip $V$ and $V_0$ 
with norms $\|\cdot\|_V = \|\cdot\| + \|\nabla \cdot\|$ 
and $\|\cdot\|_{V_0} = \|\nabla \cdot\|$, respectively. Identifying $H$ with 
its dual space $H'$ by means of the above scalar product, 
we obtain the chains of continuous 
and dense embeddings $V\subset H \subset V'$ and $V_0\subset H \subset V_0'$.
We may indicate by $\duav{\cdot,\cdot}$
the duality pairing between $V'$ and $V$, or, more generally,
between $X'$ and $X$ where $X$ is a Banach space continuously
and densely embedded into $H$. 
%
%
%
Recalling that $\bn$ stands for the outer 
unit normal vector to $\Gamma$, we also set
\begin{equation}\label{defiW}
   W := \big\{ v\in H^2(\Omega):~\dn v=0~\text{on $\Gamma$}\big\}
  \subset C^0(\overline\Omega).
\end{equation}
Then, $W$ is a closed subspace of $H^2(\Omega)$. We equip $W$ 
with the norm
\begin{equation}\label{norW}
  \| v \|_W^2 := \| v \|^2 + \| \Delta v \|^2,
\end{equation}
which (on $W$) is equivalent to the usual $H^2$-norm in view 
of well-known elliptic regularity results.

\smallskip

Next, we can fix our basic hypotheses on coefficients and data:
\beas\label{assum}
 (A1)~~$\psi \in C^2(\RR;\RR)$. 
  
  \smallskip
 
%
 
 \noindent%
 (A2)~~$g\in L^p(\Omega)$ for some $p \geq 3$.
 
  \smallskip
 
 \noindent%
 (A3)~~$z_0\in W$ with $z_0\le 1$ at every point of~$\Omega$.
 Moreover, denoting by $c\OO$ an embedding constant of 
 $W$ into $C^0(\barO)$, i.e.~a constant such that 
 $\| v \|_{C^0(\barO)}\le c\OO\| v \|_{W}$ for all $v\in W$, 
 we assume that $\epsi=\epsi(z_0):= c\OO \| 1 - z_0 \|_{W} \le 1/2$.
\eddas
It is worth commenting a bit about the above assumptions. First of all, since we 
will prove that the $z$-component of the local solution takes values in $(0,1]$,
the behavior of $\psi(r)$ for large $r$ is in fact irrelevant.
On the other hand, it may be useful to assume that
\begin{equation}\label{psi2}
 \psi(r) = r^2 \quext{for every}\, |r|\ge 2, 
\end{equation}
whence it also follows that 
\begin{equation}\label{coerc:f}
  \psi(r) \ge \frac{r^2}2 - c
   \quext{for every }\, r \in \RR.
\end{equation}
Actually, such a free ``extra-coercivity'' property will help us in the approximation
and for writing the a-priori estimates in a simpler way. 

We may also observe that~(A3) implies
\begin{equation}\label{dam:0}
  \| 1 - z_0 \|_{C^0(\barO)}
   \le c\OO \| 1 - z_0 \|_{W}
   = \epsi.
\end{equation}
Since $\epsi\le 1/2$, we have $z_0\ge 1-\varepsilon\geq 1/2$ a.e.~in~$\Omega$, i.e.~the initial
amount of damage is less than one half at (almost) any point. Of course, the ideal,
and simplest, situation occurs when $z_0\equiv 1$, i.e., the body is completely
integer at the initial time. Note that the condition $z_0\le 1$ is used only
to respect the physical significance of the model. Of
course, under such an assumption any hypothetical solution satisfies $z\le 1$ 
also for $t>0$ due to the irreversibility constraint embedded into
equation~\eqref{eq:z:strong}.

\smallskip

We can now state the main result of this paper:
\bete\label{thm:main}
 Let\/ {\rm Assumption~\ref{assum}} hold. Let $\delta\in(0,1/12]$.
 Then there exist a time $T_0\in(0,T]$ depending only on 
 $\psi$, $g$, $\epsi$ and $\delta$ and at
 least a triple $(u,z,\xi)$ of functions defined over 
 $(0,T_0)\times \Omega$ and satisfying the regularity
 properties
 \begin{align}\label{rego:u}
   & u \in C^0([0,T_0]; W^{2,\rho}(\Omega) \cap V_0) \ \mbox{ for any } \rho \in [1,p] \cap [1,6),\\
  \label{rego:z}
   & z \in H^1(0,T_0;V) \cap  C_w([0,T_0];W),\\
  \label{rego:xi}
   & \xi \in L^2(0,T_0,H), \\
  \label{sep:z}
   & c\OO \| 1 - z(t) \|_{W} \le 1 - 3\delta, \quext{for all }\, t \in [0,T_0],
 \end{align}
 where $C_w([0,T_0];X)$ stands for the space of weakly-continuous functions defined on $[0,T_0]$ with values in a Banach space $X$.
 Moreover, the triple $(u,z,\xi)$ satisfies the equations
 \begin{align}\label{eq:u}
   & -\dive(z\nabla u)=g, \\ 
  \label{eq:z}
   & \xi + z_t - \Delta z + \psi'(z) = - \frac12 | \nabla u |^2,\\
  \label{eq:xi}
   & \xi \in \alpha(z_t)      
 \end{align}
 almost everywhere in $(0,T_0)\times \Omega$, with the boundary 
 conditions~\eqref{no:flux} and the initial condition~\eqref{init}
 in the sense of traces.
In addition, if $p > 3$, then $(u,z)$ is uniquely determined by 
initial data $z_0$ and continuously depends on $z_0$. More precisely, for $i=1,2$, 
let $(u_i,z_i)$ be solutions on $[0,T_0]$. Then
$$
  \|(z_1-z_2)(t)\|_{V}+ \|(u_1-u_2)(t)\|_{V_0} \leq C \|(z_1-z_2)(0)\|_{V}
$$
for every $t \in [0,T_0]$. 
\ente
\noindent%
Note that relation \eqref{sep:z} entails in particular
\begin{equation}\label{dam:t}
  \| 1 - z(t) \|_{C^0(\barO)}
   \le c\OO \| 1 - z(t) \|_{W}
   \le 1 - 3 \delta.
\end{equation}
Hence, for any $t\in[0,T_0]$, we have $z(t,x)\ge 3 \delta >0$ for 
every $x\in \barO$. In this sense, we are able to compute a time before 
which complete damage cannot occur at any point. In such a timespan,
the system remains nondegenerate and existence of \weak solutions can
be proved. Of course, condition $\delta\le 1/12$ combined with 
assumption~(A3) implies
\begin{equation}\label{sep:16}
 c_\Omega\|1-z_0\|_{C^0(\barO)} = \epsi \le 1/2 < 3/4 \le 1-3\delta,
\end{equation}
namely there is a gap of at least $1/4$ between $1-\epsi$ and $3\delta$.
Of course the magnitude of such a gap is somehow an arbitrary
choice of ours; on the other hand, keeping it as a given value permits
us to write the estimates in a computationally simpler way.


\section{Proofs}
\label{sec:proofs}

We start with introducing a truncated version of system \eqref{eq:u:strong}-\eqref{eq:z:strong}
in the same spirit as in \cite{BS}. To this aim, for $\delta\in(0,1/12]$
we consider a mapping $T\dd\in C^{1,1}(\RR;\RR)$ such that 
\begin{equation}\label{def:Td}
  T_\delta(r) = \begin{cases}
     r & \text{~~if }r\ge 3\delta,\\
      2\delta & \text{~~if }r \le \delta
    \end{cases}
\end{equation}
and $T_\delta$ is monotone and convex in the interval $(\delta,3\delta)$  and fulfills
\begin{equation}\label{derTd}
  |T\dd'(r)|\le 1, \quad |T\dd''(r)| \le c \delta^{-1}
  \quext{for almost all }\,r\in \RR
\end{equation}
and for some $c>0$. A possible explicit choice could be 
\begin{equation}\label{def:Td2}
  T_\delta(r) = 2\delta + (4\delta)^{-1} (r-\delta)^2
  \quext{for }\, r \in(\delta,3\delta),
\end{equation}
but other options may be equally allowed. Then, the truncated system may be stated as follows:
\begin{alignat}{4}\label{eq:u:dd}
  & -\dive(T\dd(z)\nabla u)=g, &&\qquext{in }\,(0,T)\times \Omega;\\
 \label{eq:z:dd}
  & \alpha(z_t) + z_t - \Delta z + \psi'(z) \ni - \frac{T\dd'(z)}2 | \nabla u |^2,
    &&\qquext{in }\,(0,T)\times \Omega,
\end{alignat}
where, as before, the differential inclusion \eqref{eq:z:dd} may be interpreted
as the equality
\begin{equation} \label{eq:z:dd:2}
  \xi + z_t - \Delta z + \psi'(z) = - \frac{T\dd'(z)}2 | \nabla u |^2,
\end{equation}
for a suitable $\xi$ satisfying \eqref{eq:xi} at almost every point
of the parabolic cylinder.

\smallskip

We postpone to the next section a proof of the fact that a  global in time solution 
$(u,z)$ to \eqref{eq:u:dd}-\eqref{eq:z:dd} plus the initial and boundary conditions
exists in a suitable regularity class. In this part we just show
that such a solution complies with a number of a priori estimates. 
The compatibility of the estimates with the approximation will also be 
discussed later on.
In this procedure, we will denote by $c$ a generic
positive constant depending only on the assigned data of the problem,
including $\epsi$ and the final time $T$. On the other hand, $c$ will
{\sl not}\/ be allowed to depend on $\delta$ (so when $\delta$ appears
in the computations, it will be kept explicit).

Our purpose is to construct in a computable way a time interval
$(0,T_0)$, with $T_0>0$ possibly smaller than $T$ and 
depending on the given constants $\delta$ and $\epsi$,
such that $u(t,x)\ge 3 \delta$ for a.e.~$(t,x)\in(0,T_0)\times \Omega$. In
this way, due to \eqref{def:Td}, $(u,z)$ will turn out to solve
the original system \eqref{eq:u:strong}-\eqref{eq:z:strong} 
in that time span.

To start, we perform the analogue of the energy estimate
described in the introduction. Testing \eqref{eq:u:dd} by $u_t$,
\eqref{eq:z:dd} by $z_t$, and performing standard manipulations 
(note in particular that the product $\xi z_t$
is a.e.~equal to $0$ since $\alpha(z_t)$ may contain nonzero values 
only at $z_t=0$), we easily arrive at
\begin{equation}\label{eq:E:dd}
  \ddt \calE\dd(t) + \| z_t \|^2 = 0,
\end{equation}
with the truncated energy functional
\begin{equation}\label{Edd}
   \calE\dd(t) = \io \Big( \frac{T\dd(z)}2 | \nabla u |^2 - gu
    + \frac{1}2 |\nabla z |^2 + \psi(z) \Big).
\end{equation}
Note that the use of test functions and the integrations by parts 
performed to deduce this estimate and the subsequent ones will be justified
as far as one works with the regularized solutions (see the next 
section for details).
Now, as we integrate \eqref{eq:E:dd} over some time interval $(0,t)$, we see
that $\calE\dd(0)$ also depends on the ``initial value''
$u_0=u|_{t=0}$. However, in view of the quasi-static nature of the system,
$u_0$ is not a datum, but has to be computed by evaluating \eqref{eq:u:dd} at the time $t=0$.
Namely, $u_0$ corresponds to the (unique)
solution to the elliptic problem
\begin{equation}\label{co:09}
  -\dive(T\dd(z_0)\nabla u_0)=g, \qquext{in }\,\Omega,
\end{equation}
complemented with the  homogeneous Dirichlet boundary condition. 
In view of Assumption~(A3) and of the
fact $3 \delta \le 1-\epsi$, we actually have
$T\dd(z_0)=z_0 \ge 1/2$. Hence, testing \eqref{co:09}
by $u_0$, we obtain
\begin{equation}\label{co:10}
  \frac12 \| \nabla u_0 \|^2
   \le \io T\dd(z_0) |\nabla u_0| ^2 
    = (g,u_0)
   \le \| g \| \| u_0 \|
   \le \frac14 \| \nabla u_0 \|^2 + c,
\end{equation}
where Poincar\'e's inequality has also been used. This fact implies in 
particular that 
\begin{align}\no
  \big| \calE\dd(0) \big|
  & = \bigg| \io \Big( \frac{T\dd(z_0)}2 | \nabla u_0 |^2 - gu_0
    + \frac{1}2 |\nabla z_0 |^2 + \psi(z_0) \Big) \bigg|\\
 \label{Edd0}
  & = \bigg| \io \Big( - \frac{T\dd(z_0)}2 | \nabla u_0 |^2 
    + \frac{1}2 |\nabla z_0 |^2 + \psi(z_0) \Big) \bigg|
    \le c ( 1 + \| z_0 \|_V^2 ),
\end{align}
with $c$ independent of $\delta$. Hence, recalling that $z_0\in V$,
$z_0\le 1$ almost everywhere (cf.~Assumption~(A3)),
we see in particular that our assumptions on the initial
data imply the finiteness of the energy at $t=0$.

Integrating \eqref{eq:E:dd} over the generic time interval
$(0,t)$ (where the choice of the admissible ``small'' time $t>0$ 
will be made clear later on), we then infer  that
\begin{equation}\label{eq:E:dd:2}
  \calE\dd(t) + \itt \| z_t \|^2 = 
   \calE\dd(0) \le c ( 1 + \| z_0 \|_V^2 ).
\end{equation}
Now, using Poincar\'e's inequality, we arrive at
\begin{equation}\label{co:12}
  \bigg| \io gu \bigg| 
   \le \| g \| \| u \| 
   \le c \| g \| \| \nabla u \| 
   \le \frac\delta2 \| \nabla u \|^2
   + \frac{c}{\delta}.
\end{equation}
As a consequence of the above relations \eqref{coerc:f} and \eqref{def:Td}, we have
\begin{equation}\label{eq:E:dd:3}
  \calE\dd(t) \ge \frac\delta2 \| \nabla u(t) \|^2
   + \frac12 \| z(t) \|_V^2 
   - \frac{c}{\delta}.
\end{equation}
Combining \eqref{eq:E:dd:2} with \eqref{eq:E:dd:3}, we then obtain
the a priori estimates,
\begin{align}\label{st:11}
  & \| z \|_{L^\infty(0,t;V)} \le c\big( \delta^{-1/2} + \| z_0 \|_V \big),\\
 \label{st:12}
  & \| u \|_{L^\infty(0,t;V_0)} \le c \delta^{-1/2} \big( \delta^{-1/2} + \| z_0 \|_V \big),\\
 \label{st:13}
  & \| z_t \|_{L^2(0,t;H)} \le c \big( \delta^{-1/2} + \| z_0 \|_V \big).
\end{align}
Next, evaluating \eqref{eq:u:dd} at the generic time $t$ and testing
it by $u$, applying once more Poincar\'e's inequality, we obtain
\begin{align}\no
  \| \nabla u \|^2
      & = \io |\nabla u| ^2  = \io \frac{T\dd(z)}{T\dd(z)}|\nabla u|^2 
    \le  \left\| \frac1{T\dd(z)} \right\|_{L^\infty(\Omega)}
       \io T\dd(z) |\nabla u|^2\\
 \label{co:20}    
   &    = \left\| \frac1{T\dd(z)} \right\|_{L^\infty(\Omega)}
       (g,u)
 \le c \left\| \frac1{T\dd(z)} \right\|_{L^\infty(\Omega)}
       \| g \| \| \nabla u \|,
\end{align}
whence 
\begin{equation}\label{co:14}
  \| \nabla u \|
   \le c \left\| \frac1{T\dd(z)} \right\|_{L^\infty(\Omega)},
\end{equation}
with computable $c>0$ also depending on $g$.

Now let us define, for $r\in \RR$, 
\begin{equation}\label{co:21}
  \phi\dd(r):= \frac{1}{T\dd(1-r)},
  \qquext{so that }\,
   \frac{1}{T\dd(r)} 
   = \frac{1}{T\dd(1-(1-r))}
   = \phi\dd(1-r).
\end{equation}
In other words,  for $r\in \R$, the function $\phi\dd(r)$ is a regularization of the function $r \mapsto 1/(1-r)_+$; 
in particular, $\phi_\delta(r)=(1-r)^{-1}$ for $r \leq 1-3\delta$. Notice also that $\phi\dd$ is 
non-decreasing on  $\R$. 

 By the use of \eqref{co:21}, \eqref{co:14} can be rewritten as 
\begin{equation}\label{co:22}
  \| \nabla u \|
   \le c \| \phi\dd(1-z) \|_{L^\infty(\Omega)}
   =  c \phi\dd\big( \|1-z \|_{L^\infty(\Omega)} \big).
\end{equation}
Next, let us observe that \eqref{eq:u:dd} may be equivalently rewritten
as 
\begin{equation}\label{eq:u:dd:3}
  - T\dd(z) \Delta u = g + T\dd'(z) \nabla z \cdot \nabla u.
\end{equation}
We now compute the $L^2$- and $L^3$-norms of both sides
of the above relation. Observing that $T\dd(r)\ge 2\delta$ with 
$|T\dd'(r)| \le 1$ for every $\delta\in(0,1/12]$
and $r\in\RR$, and using elementary interpolation and embedding inequalities  along with \eqref{norW},
we first find that 
\begin{align}\no
  2 \delta \| \Delta u \|
   & \le \| g \| + \| \nabla z \|_{L^6(\Omega)} \| \nabla u \|_{L^3(\Omega)} \\
 \no
   & = \| g \| + \| \nabla ( z - 1) \|_{L^6(\Omega)} \| \nabla u \|_{L^3(\Omega)} \\
 \no   
   & \le \| g \| + c \| z - 1 \|_W \| \nabla u \|^{1/2} \| \Delta u \|^{1/2} \\   
 \no   
   & \le c + c \delta^{-1/2} \big( \| z - 1 \| + \| \Delta z \| \big) \| \nabla u \|^{1/2} \delta^{1/2} \| \Delta u \|^{1/2} \\   
 \label{new:23}
  & \le c + c \delta^{-1} \big( \| z - 1 \|^2 + \| \Delta z \|^2 \big) \| \nabla u \| 
   + \delta \| \Delta u \|.
\end{align}
Analogously, combining the Gagliardo-Nirenberg inequality \cite{Nir} with standard
elliptic regularity results of $L^p$-type, we infer that
\begin{equation}\label{gani}
  \| \nabla v \|_{L^6(\Omega) }
   \le c \| \Delta v \|_{L^3(\Omega)}^{2/3} \| \nabla v \|^{1/3}, 
\end{equation}
which holds for every $v\in V_0 \cap W^{2,3}(\Omega)$. Using such a relation, we deduce that
\begin{align}\no
  2 \delta \| \Delta u \|_{L^3(\Omega)}
   & \le \| g \|_{L^3(\Omega)} + \| \nabla z \|_{L^6(\Omega)} \| \nabla u \|_{L^6(\Omega)} \\
 \no
   & = \| g \|_{L^3(\Omega)} + \| \nabla ( z - 1) \|_{L^6(\Omega)} \| \nabla u \|^{1/3} \| \Delta u \|_{L^3(\Omega)}^{2/3}\\
 \no   
   & \le c + c \delta^{-2/3} \big( \| z - 1 \| + \| \Delta z \| \big)  \| \nabla u \|^{1/3} \delta^{2/3} \| \Delta u \|_{L^3(\Omega)}^{2/3}\\
 \label{new:23b}
  & \le c + c \delta^{-2} \big( \| z - 1 \|^3 + \| \Delta z \|^3 \big) \| \nabla u \| 
   + \delta \| \Delta u \|_{L^3(\Omega)}.
\end{align}
Hence, recalling also \eqref{co:22}, \eqref{new:23} and \eqref{new:23b} imply respectively
\begin{align}\label{sti:l2}
  & \| \Delta u \| \le c \delta^{-1} + c \delta^{-2} \big( \| z - 1 \|^2 + \| \Delta z \|^2 \big) 
      \phi\dd\big( \|1-z \|_{L^\infty(\Omega)} \big),\\ 
 \label{sti:l3}
  & \| \Delta u \|_{L^3(\Omega)} \le c \delta^{-1} + c \delta^{-3} \big( \| z - 1 \|^3 + \| \Delta z \|^3 \big) 
      \phi\dd\big( \|1-z \|_{L^\infty(\Omega)} \big).
\end{align}
As a next step, we test \eqref{eq:z:dd} by $-\Delta z_t$. Then, using the monotonicity of $\alpha$ 
and the no-flux boundary conditions, we would expect that
\begin{equation}\label{sep:21}
   (\alpha(z_t),-\Delta z_t) 
    = \io \alpha'(z_t) | \nabla z_t |^2 
    \ge 0.
\end{equation}
On the other hand, the above computation is formal. Indeed, $\alpha$ is a nonsmooth maximal
monotone graph (and $\alpha(z_t)$ has to be interpreted as a selection $\xi$
(cf.~\eqref{eq:xi}). Nevertheless, the inequality $(\xi,-\Delta z_t) \ge 0$ is 
valid anyway, and it could be rigorously proved by proceeding, e.g., along the lines of~\cite[Lemma~2.4]{SP}
(see also Remark~\ref{zt} below for a further justification of this procedure). 
Hence, we deduce  that
\begin{equation}\label{co:26}
  \| \nabla z_t \|^2 + \frac12 \ddt \| \Delta z \|^2
    \le \io | \psi''(z) \nabla z \cdot \nabla z_t | 
     + \frac12 \io \big| \nabla \big( T\dd'(z) | \nabla u | ^2 \big) \cdot \nabla z_t \big| 
     =: I_1 + I_2
\end{equation}
and we need to control the terms on the \rhs. First of all,
by~(A1) and \eqref{psi2} we have that $\psi''\in L^\infty(\RR)$,
whence
\begin{equation}\label{co:27}
  I_1 \le c \| \nabla z \|  \| \nabla z_t \| 
   \le \frac 16 \| \nabla z_t \|^2 + c \| 1 - z \|_V^2.
\end{equation}
Next, recalling \eqref{derTd}, we easily obtain
\begin{equation}\label{co:28}
   I_2 \le  c \delta^{-1} \io | \nabla u |^2 |\nabla z \cdot \nabla z_t |
   + c \io | D^2 u | | \nabla u | | \nabla z_t | =: I_{2,1} + I_{2,2}.
\end{equation}
Furthermore, using also \eqref{sti:l2}, we infer that
\begin{align}\no
  I_{2,1} & \le  c \delta^{-1} \| \nabla u \|_{L^6(\Omega)}^2
   \| \nabla z \|_{L^6(\Omega)} \| \nabla z_t \|\\
 \no
   & \le  c \delta^{-1} \| \Delta u \|^2
    \big ( \| z-1 \| + \| \Delta z \| \big) \| \nabla z_t \|\\
 \no
   & \le \frac16 \| \nabla z_t \|^2 + 
    c \delta^{-2} \| \Delta u \|^4
    \big ( \| z-1 \|^2 + \| \Delta z \|^2 \big)\\
 \no
   & \le \frac16 \| \nabla z_t \|^2 + c \delta^{-6} \big ( \| z-1 \|^{2} + \| \Delta z \|^{2} \big)\\
 \label{co:28b}  
  & \mbox{}~~~~~~~~~~~~~~~~~~~~~~
 + c \delta^{-10} \big( \| z - 1 \|^{10} + \| \Delta z \|^{10} \big) 
         \phi\dd^4\big( \|1-z \|_{L^\infty(\Omega)} \big). 
\end{align}
Similarly, using elliptic regularity  along with \eqref{sti:l2} and \eqref{sti:l3} as well as Young's inequality, we obtain
\begin{align}\no
  I_{2,2} & \le c \| D^2 u \|_{L^3(\Omega)}
   \| \nabla u \|_{L^6(\Omega)} \| \nabla z_t \|\\
 &
   \le \frac16 \| \nabla z_t \|^2 + c \| \Delta u \|^2_{L^3(\Omega)} \| \Delta u \|^2 
\nonumber\\
   & \le \frac16 \| \nabla z_t \|^2 + 
    c \delta^{-10} \big ( \| z-1 \|^{10} + \| \Delta z \|^{10} \big)
     \phi\dd^4\big( \|1-z \|_{L^\infty(\Omega)} \big)
 \nonumber \\
 & \quad + c \delta^{-6} \big ( \| z-1 \|^4 + \| \Delta z \|^4 \big)
     \phi\dd^2\big( \|1-z \|_{L^\infty(\Omega)} \big) \nonumber \\
 &\quad + c \delta^{-8} \big ( \| z-1 \|^6 + \| \Delta z \|^6 \big)
     \phi\dd^2\big( \|1-z \|_{L^\infty(\Omega)} \big)
    + c \delta^{-4}\nonumber \\
 & \leq \frac16 \| \nabla z_t \|^2 + 
    c \delta^{-10} \big ( 1+ \| z-1 \|_W^{10} \big)
     \left[1+ \phi\dd^4\big( \|1-z \|_{L^\infty(\Omega)} \big) \right]
  \label{co:29}  
\end{align}
for $\delta \in (0,1/12]$. 
Notice that this is actually the only point 
in the existence proof where we need the control on the $L^3$-norm
of $\Delta u$ (and, in turn, the assumption $g\in L^3(\Omega)$).

Collecting \eqref{co:26}-\eqref{co:29} gives
\begin{align}\no
  \| \nabla z_t \|^2 
   + \ddt \| \Delta z \|^2
    & \le c \| 1 - z \|_V^2
     + c \delta^{-6} \| z-1 \|_W^2 \\
 \label{co:30}             
    & \quad + c \delta^{-10} \left( 1+ \| z-1 \|_W^{10} \right)
 \left[ 1+ \phi\dd^4\big( \| 1-z \|_{L^\infty(\Omega)} \big) \right].
\end{align}
In order to deduce some useful information from the above relation, we observe 
the inequality
\begin{equation}\label{co:31}
   \ddt \| 1 - z \|^2  
   \le 2 |(1-z,z_t)| \le c \| 1 - z \|^4  + c \| z_t \|^{4/3}.
\end{equation}
Adding it to \eqref{co:30} and rearranging terms, with the aid of Young's inequality, we arrive at
\begin{align}
   \ddt \| 1 - z \|_W^2
   + \| \nabla z_t \|^2 
  \le 
 \label{co:31x}
    c \delta^{-10} \left( 1+\| z-1 \|_W^{10} \right)
 \left[ 1+\phi\dd^4\big( \| 1-z \|_{L^\infty(\Omega)} \big) \right] 
    + c \| z_t \|^{4/3}.
\end{align}
Let us now multiply the above by $c\OO^2$, the embedding constant
of $H^2(\Omega)$ into $C^0(\barO)$ as introduced before. Then, setting
\begin{equation}\label{def:y}
  y(t):= c\OO^2 \| 1 - z(t) \|_W^2 \stackrel{\eqref{dam:t}}\geq \|1-z(t)\|_{L^\infty(\Omega)}^2, 
\end{equation}
and temporarily neglecting the nonnegative term $\| \nabla z_t \|^2$
on the \lhs, 
we deduce the differential inequality
\begin{align}\no
  y'(t) 
\label{co:33}    
   & \le c_1 \delta^{-10} \left[ 1 + y^5(t) \right] \left[ 1 + \phi\dd^4( y^{1/2}(t)) \right]
    + c_2 \| z_t \|^{4/3},
\end{align}
where it is worth noting that $y_0:=y(0) = \epsi^2 \le 1/4$ by assumption~(A3) and 
$c_1$, $c_2$ are computable positive constants independent of $\delta$.
Dividing both sides by $[ 1 + y^5(t) ] [ 1 + \phi\dd^4( y^{1/2}(t)) ]$,
which is clearly larger than~$1$, we then obtain
\begin{equation}\label{co:34}
  \ddt B\dd(y) :=  \frac{1}{\left[ 1 + y^5(t) \right] \big[ 1 + \phi\dd^4( y^{1/2}(t)) \big]} y'
   \le c_1 \delta^{-10} + c_2 \| z_t \|^{4/3},
\end{equation}
where the function $B\dd$ is defined by the \lhs, namely we have set
\begin{equation}\label{co:34b}
  B\dd(s) :=  \int_0^s \frac{\dir}{( 1 + r^5 ) \big[ 1 + \phi\dd^4( r^{1/2} ) \big]}.
\end{equation}
Here, we note that 
the function $B\dd$, as far as $\delta$ 
is a fixed number in the given range $(0,1/12]$, is well defined and strictly increasing on $\R$. 
Now, it is clear that, for $s\in[0,1]$,
\begin{equation}\label{sep:11}
   \frac{1}2 \int_0^s \frac{\dir}{1 + \phi\dd^4( r^{1/2} )}  
    \le B\dd(s) \le \int_0^s \frac{\dir}{1 + \phi\dd^4( r^{1/2} )}.
\end{equation}
Moreover, from \eqref{co:21}, we observe that, for $r^{1/2}\in [0,1-3\delta]$,
or equivalently $r\in [0,(1-3\delta)^2]$,
\begin{equation}\label{sep:12}
  \frac{1}{1 + \phi\dd^4( r^{1/2} )} 
   = \frac{(1-r^{1/2})^4}{(1-r^{1/2})^4+1},
\end{equation}
whence, we can notice that, as far as $s$ lies in the range $[0,(1-3\delta)^2]$, 
the expression of $B\dd(s)$ is {\bf independent of $\delta$}
so that, for such $s$, we can simply write $B(s)$ in place of $B\dd(s)$.
Notice also that, at largest, $\delta=1/12$; hence
$(1-3\delta)^2$ is always at least $9/16$. 

Integrating \eqref{co:34} in time and using \eqref{st:13} with H\"older's inequality,
we obtain
\begin{align}\no
  B\dd(y(t)) & \le B\dd(y_0) + \int_0^t \big( c_1 \delta^{-10} + c_2 \| z_t \|^{4/3} \big)
   \le B\dd(\epsi^2) + c_1 \delta^{-10} t + c_2 t^{1/3} \big( \delta^{-1} + \| z_0 \|_V^2 \big)^{2/3}\\
 \label{co:35n}
   & \le B\dd(\epsi^2) + c_3 \delta^{-10} t^{1/3},
\end{align}
where the new constant $c_3$ may also depend on $z_0$ and $T$.

On the other hand, due to \eqref{sep:16} along with the strict increase of $B_\delta$, 
\eqref{co:35n} can be rewritten as
\begin{equation}\label{sep:15}
  y(t) \le B\dd^{-1} \big( B\dd(\epsi^2) + c_3 \delta^{-10} t^{1/3} \big)
   = B\dd^{-1} \big( B(\epsi^2) + c_3 \delta^{-10} t^{1/3} \big),
\end{equation}
where we used that $\epsi^2\le 1/4 < 9/16\le (1-3\delta)^2$.

Now, since $\delta$ is assigned and $c_3$ is a computable constant
depending only on the given parameters of the system, using that
$B\dd$ is strictly monotone (hence such is its inverse $B\dd^{-1}$),
we deduce that there exists 
$T_0\in(0,T]$ so small that, for every $t\in[0,T_0]$, there holds
\begin{equation}\label{sep:15b}
  B(\epsi^2) + c_3 \delta^{-10} t^{1/3} 
   \le B((1-3\delta)^2) = B\dd((1-3\delta)^2).
\end{equation}
In other words $T_0$ can be defined as the largest time $t \in (0,T] $ such that
$B(\epsi^2) + c_3 \delta^{-10} t^{1/3} \leq B((1-3\delta)^2)$, that is, 
$$
T_0 = \left( \frac{B((1-3\delta)^2)-B(\varepsilon^2)}{c_3 \delta^{-10}}\right)^3 \wedge T \in (0,T].
$$
As a consequence, in the range $[0,T_0]$ the expression of $B\dd$ is independent of $\delta$ and 
\eqref{sep:15} reduces to 
\begin{equation}\label{sep:15d}
  y(t) \le B^{-1} \big( B(\epsi^2) + c_3 \delta^{-10} t^{1/3} \big)
   \le  B^{-1} \big( B((1-3\delta)^2) \big), \quext{for all }\,t\in[0,T_0],
\end{equation}
which in turn implies 
\begin{equation}\label{co:40}
  \| 1 - z(t) \|_{C^0(\barO)} 
   \le c\OO \| 1 - z(t) \|_W 
   = y^{1/2}(t)
   \le 1 - 3\delta
\end{equation}
and consequently
\begin{equation}\label{co:41}
  z(t,x) \ge 3 \delta \quext{for all~}\,t\in [0,T_0],~x\in \overline{\Omega}.
\end{equation}
This entails in particular that, for every $t\in [0,T_0]$, there holds
$T\dd(z(t)) = z(t)$ a.e.~in $\Omega$, whence $(u,z)$ turns out to solve
the original system \eqref{eq:u:strong}-\eqref{eq:z:strong}.

%
%
%
%
%
%

We finally prove the regularity properties \eqref{rego:u}-\eqref{rego:xi}.
First of all, we shall check \eqref{rego:z}; the fact $z \in C_w([0,T_0];W)$ 
comes from \eqref{co:40}, while $z \in H^1(0,T_0;V)$ follows from integration of \eqref{co:31x} over $(0,T_0)$. 

Next, we prove \eqref{rego:u} which is a bit more tricky.
First of all, let $(z_i,u_i)$ be two solutions 
for \eqref{eq:u:strong}, \eqref{eq:z:strong} on $[0,T_0]$. Then by subtraction, we have
$$
  -\dive \left[ z_1(\nabla u_1-\nabla u_2) + (z_1-z_2)\nabla u_2 \right]=0,
    \quext{in }\, \Omega.
$$
Test it by $u_1-u_2$. We see that
\begin{align*}
 \int_\Omega z_1 |\nabla (u_1-u_2)|^2 
&=  - \int_\Omega(z_1-z_2) \nabla u_2 \cdot \nabla (u_1-u_2)\\
&\leq \|z_1-z_2\|_{L^4(\Omega)} \|\nabla u_2\|_{L^4(\Omega)} \|\nabla (u_1-u_2)\|,
\end{align*}
which entails
$$
  3\delta \|\nabla (u_1-u_2)\| \leq \|z_1-z_2\|_{L^4(\Omega)} \|\nabla u_2\|_{L^4(\Omega)}.
$$
Hence we may conclude in particular that
\begin{equation}\label{el-subt}
  3 \delta \|\nabla u(t) - \nabla u(s)\| \leq \|z(t)-z(s)\|_{L^4(\Omega)} 
   \sup_{\tau \in [0,T_0]} \|u(\tau)\|_{H^2(\Omega)}, \quext{for }\, t,s \in [0,T_0],
\end{equation}
and, therefore, $t \mapsto u(t)$ turns out to be continuous on $[0,T_0]$ with values in $V_0$.
Furthermore, \eqref{eq:u:strong} implies
\begin{equation}\label{eq:u:div}
  - \Delta u = \dfrac g z + \frac{\nabla z}z \cdot \nabla u \ \mbox{ in } (0,T_0) \times \Omega.
\end{equation}
Note that $t \mapsto 1/z(t)$ is continuous with values in $L^\infty(\Omega)$ on $[0,T_0]$ 
(indeed, $H^1(0,T_0;V) \cap L^\infty(0,T_0;W)$ is embedded in $C^0([0,T_0];L^\infty(\Omega))$ 
and $z$ is uniformly away from zero in $(0,T_0) \times \Omega$). Since $u(t)$ is also bounded 
in $W^{2,3}(\Omega)$ for any $t \in [0,T_0]$ and $u \in C^0([0,T_0];V_0)$, the map $t \mapsto \nabla u(t)$ 
is continuous on $[0,T_0]$ strongly in $L^q(\Omega)$ for any $q \in [1,+\infty)$. 
On the other hand, thanks to an Aubin-Lions type embedding (see, e.g.,~\cite{Si}), we may observe that
$$
  L^\infty(0,T_0;H^1(\Omega)) \cap H^1(0,T_0;H) \hookrightarrow C^0([0,T_0];L^q(\Omega)),
   \quext{for any }\,q\in [1,6).
$$  
Applying this to $\nabla z$, we can verify that $t \mapsto \nabla z(t)$ is of class 
$C^0([0,T_0];L^q(\Omega))$ for $q\in [1,6)$. Combining the above facts, we deduce that
$t \mapsto z^{-1}(t) \nabla z(t) \cdot \nabla u(t)$ is continuous strongly 
in $L^q(\Omega)$ for any $q \in [1, 6)$. Thus the (strong) continuity of 
$t \mapsto \Delta u(t)$ in $L^\rho(\Omega)$ for any $\rho \in [1,p] \cap [1,6)$ 
on $[0,T_0]$ follows from (A2), \eqref{eq:u:div} and the facts observed so far. 

Concerning the continuous dependence of solutions on the initial data, 
let $(u_i,z_i)$ for $i=1,2$ be two solutions on $[0,T_0]$ and assume~(A2)
holds for $p>3$. Then, setting 
$Z = z_1-z_2$ and $U = u_1-u_2$, by subtraction, we have
$$
  \alpha(\partial_t z_1) - \alpha(\partial_t z_2) + Z_t - \Delta Z + \psi'(z_1) - \psi'(z_2) 
   \ni - \frac 1 2 \left( |\nabla u_1|^2 - |\nabla u_2|^2 \right). 
$$
Test both sides by $Z_t$ and employ the monotonicity of $\alpha$. Moreover, note that $u_i \in L^\infty(0,T_0;W^{2,\rho}(\Omega))$, $i=1,2$,
where now $\rho > 3$, and the embedding $W^{1,\rho}(\Omega) \hookrightarrow L^\infty(\Omega)$. We then obtain
\begin{equation}\label{co:81}
  \frac 1 2 \|Z_t\|^2 + \frac 1 2 \ddt \|\nabla Z\|^2 \leq \|\psi'(z_1)-\psi'(z_2)\|^2 
   + \frac 1 4 \|\nabla (u_1+u_2)\|_{L^\infty(\Omega)}^2 \|\nabla U\|^2
\end{equation}
by using
\begin{align*}
  \left| \io \left( |\nabla u_1|^2 - |\nabla u_2|^2 \right) Z_t \right|
   &\leq \frac 1 2 \|Z_t\|^2 + \frac 1 2 \left\| (\nabla u_1 + \nabla u_2) \cdot \nabla U \right\|^2\\
   &\leq \frac 1 2 \|Z_t\|^2 + \frac 1 2 \left\| \nabla u_1+\nabla u_2 \right\|_{L^\infty(\Omega)}^2 \|\nabla U\|^2.
\end{align*}
Next, notice that 
$$
  \|\psi'(z_1)-\psi'(z_2)\| \leq c \|Z\|
$$
for some constant $c>0$. Hence, \eqref{co:81} implies
\begin{align*}
  \frac 1 2 \|Z_t\|^2 + \frac 1 2 \ddt \|\nabla Z\|^2 \leq c \left( \|Z\|^2 + \|\nabla U\|^2\right),
\end{align*}
which, along with \eqref{el-subt}, implies
\begin{align*}
  \frac 1 2 \|Z_t\|^2 + \frac 1 2 \ddt \|\nabla Z\|^2 \leq 
   c \left( \|Z\|^2 + \|Z\|_{L^4(\Omega)}^2\right) \leq c \|Z\|_{V}^2.
\end{align*}
Summing the elementary inequality
\begin{equation}\label{co:82}
  \frac12 \ddt \| Z \|^2 
   \le \frac14 \| Z_t \|^2 + \| Z \|^2    
\end{equation}
in order to recover the full $V$-norm on the \lhs\ and subsequently 
using Gronwall's lemma, we conclude that
$$
  \|Z(t)\|_{V}^2 \leq c \|Z(0)\|_{V}^2 \quext{for }\, t \in [0,T_0].
$$
Moreover, \eqref{el-subt} yields
$$
  \|U(t)\|_{V_0}^2 \leq C \|Z(t)\|_{V}^2\ \quext{for }\, t \in [0,T_0].
$$
The uniqueness follows immediately under the assumption $Z(0)=0$, i.e., when the initial
data are the same.

Finally, \eqref{rego:xi} follows from \eqref{rego:u}-\eqref{rego:z} 
and a comparison of terms in \eqref{eq:z}.
This concludes the proof of Theorem~\ref{thm:main} provided that we
can exhibit a regularization of the system for which:
\begin{itemize}
\item we can prove existence of sufficiently smooth solutions
on the time interval $(0,T)$;
\item we can show compatibility of the regularization 
with the a priori estimates performed above.
\end{itemize}
This will be the purpose of the next section.


\section{Approximation}
\label{sec:appro}

We introduce here a regularization of system \eqref{eq:u}-\eqref{eq:xi}
for which existence 
can be proved by means of a fixed point argument. Namely, letting 
$\epsilon\in (0,1)$ be a regularization parameter intended to go
to $0$ in the limit, we introduce the system
\begin{alignat}{4}\label{eq:u:ee}
& \epsilon \Delta^2 u -\dive(T\dd(z)\nabla u)=g, && \qquext{in }\,(0,T)\times \Omega,\\
 \label{eq:z:ee}
& \alpha(z_t) + z_t - \Delta z + \psi'(z) \ni - \frac{T\dd'(z)}2 | \nabla u |^2,&& \qquext{in }\,(0,T)\times \Omega,
\end{alignat}
(for brevity, here we avoid to introduce the notation $\xi$ for the
representative of $\alpha(z_t)$, cf.~\eqref{eq:xi}). It is worth observing that, in this approximation,
we {\sl do not need}\/ to smooth out the operator $\alpha$. Hence, 
the irreversibility constraint and the related property will hold also 
for solutions to \eqref{eq:u:ee}-\eqref{eq:z:ee}.

The above relations are complemented with the same initial and boundary condition
considered before and with the additional boundary condition
\begin{equation}\label{co:51}
  \Delta u = 0, \quext{on }\,(0,T)\times \Gamma.
\end{equation}
It is worth noting from the very beginning that the system above is fully compatible
with the local a-priori estimates performed in the previous section. Indeed, 
as we test \eqref{eq:u:ee} by $u_t$ we obtain an additional (positive) term in
the energy functional, namely we have
\begin{equation}\label{Eeedd}
   \calE_{\epsilon,\delta}(t) = \io \Big( \frac\epsilon2|\Delta u |^2 
    + \frac{T\dd(z)}2 | \nabla u |^2 - gu
    + \frac{1}2 |\nabla z |^2 + \psi(z) \Big),
\end{equation}
and the new term is a source of additional {\sl a-priori}\ regularity.
On the other hand, the elliptic regularization is also compatible with the procedure
used to get the differential inequality \eqref{co:34}. Actually, the key
estimates \eqref{new:23} and \eqref{new:23b} can still be obtained 
similarly as before. Namely, to get the analogue of \eqref{new:23} we 
now need to test \eqref{eq:u:ee} by $-\Delta u$, whereas 
for \eqref{new:23b} we test \eqref{eq:u:ee} by $-|\Delta u|\Delta u$
and notice that
\begin{equation}\label{co:52}
  \io - \epsilon \Delta^2 u ( | \Delta u | \Delta u )
   = 2 \epsilon \io | \Delta u | |\nabla \Delta u |^2 \ge 0,
\end{equation}
also in view of the additional boundary condition \eqref{co:51}.

On the other hand, the new term provides additional compactness and it may help
to solve \eqref{eq:u:ee}-\eqref{eq:z:ee} by means of a fixed point argument.
We now sketch a possible procedure (which, in some sense, is inspired by
the argument given in \cite{BS}), leaving the details to the reader.

\smallskip

{\bf (1)~~}We take a prescribed function $\baru$ instead of $u$ in \eqref{eq:z:ee}. More precisely,
we choose 
\begin{equation}\label{reg:ubar}
 \baru \in L^4(0,T;W^{2,3}(\Omega) \cap V_0).
\end{equation}
This in particular implies that 
\begin{equation}\label{reg:ubar2}
  | \nabla \baru |^2 \in L^2(0,T;V)
\end{equation}
as a direct check shows. The corresponding equation 
\begin{equation}\label{eq:fp2}
  \alpha(z_t) + z_t - \Delta z + \psi'(z) \ni - \frac{T\dd'(z)}2 | \nabla \baru |^2
\end{equation}
is a parabolic equation with the Lipschitz nonlinearity $T\dd'(z)$ and the
nonsmooth term $\alpha(z_t)$. For this type of equation the regularity theory 
is well-established. For instance, one can test it by $-\Delta z_t$
(see also Remark~\ref{zt} below). 
Then, using the monotonicity of $\alpha$, condition \eqref{reg:ubar}, the 
Lipschitz continuity of $T\dd'$, and Gronwall's lemma, one may
deduce the existence of at least one solution $z$ in the same 
regularity class of Theorem~\ref{thm:main}, namely
\begin{equation}\label{reg:z:fp}
  z \in H^1(0,T;V) \cap L^\infty(0,T;W).
\end{equation}
Moreover, such a solution is readily seen to be unique. To check this fact
it suffices to take a couple of solutions (with the same proposed 
$\baru$), compute correspondingly the difference
of \eqref{eq:fp2}, and test it by the difference of the $z_t$'s.
Then, exploiting the monotonicity of $\alpha$ one can easily obtain a 
contraction estimate.

\smallskip

{\bf (2)~~}We plug the function $z$ obtained at the previous step into \eqref{eq:u:ee}. 
This gives rise to a fourth order elliptic equation, whose leading term is linear, 
with the boundary conditions $u=\Delta u = 0$ on $(0,T) \times \Gamma$. Hence, it has a unique 
weak solution $u \in L^\infty(0,T;H^2(\Omega) \cap V_0)$.
Moreover, we can also prove that
\begin{equation}\label{reg:u:fp}
u \in L^\infty(0,T;H^{ 4 }(\Omega)).
\end{equation}
Indeed, rewrite \eqref{eq:u:ee} as
\begin{equation}\label{v-ee} 
  \epsilon \Delta^2 u - T_\delta(z) \Delta u = T_\delta'(z)\nabla z \cdot \nabla u + g \ \mbox{ in } (0,T) \times \Omega, 
\end{equation}
which is complemented with the homogeneous Dirichlet boundary conditions and where the right-hand side lies at least on $L^\infty(0,T;L^2(\Omega))$. Hence the $L^2$-regularity theory for higher-order elliptic operators entails $u(\cdot,t) \in H^4(\Omega)$ for a.e.~$t \in (0,T)$. More precisely, we can set $v = - \Delta u$ and apply the $L^2$ elliptic regularity of second order type. Then we have
$$
\esssup_{t \in (0,T)} \int_\Omega |\partial_{ij}^2 v|^2 \leq C,
$$
where $\partial_{ij} = \partial^2/\partial x_i \partial x_j$ for $i,j = 1,2,3$. Here we used $u \in L^\infty(0,T;V_0)$ and \eqref{reg:z:fp} along with $W \subset L^\infty(\Omega)$. Using relation $v = -\Delta u$ and integrating by parts, the above can be rewritten as
$$
\esssup_{t \in (0,T)}\int_\Omega |\partial_{ijkl}^4 u|^2 \leq C \quad \mbox{ for }\ i,j,k,l=1,2,3,
$$
which yields $u \in L^\infty(0,T;H^4(\Omega))$.

\smallskip

{\bf (3)~~}We finally consider the mapping $\baru \mapsto u$ and we aim to apply the
Schauder fixed point theorem to this map in order to get existence of at least one 
local in time solution to the initial-boundary value problem for 
\eqref{eq:u:ee}-\eqref{eq:z:ee}. The most delicate point is proving 
compactness, because the system is quasi-stationary and we have no
information on $u_t$. On the other hand, by \eqref{reg:z:fp} and the Aubin-Lions
theorem, one can easily obtain that the mapping 
$\baru \mapsto z$ is completely continuous from the space \eqref{reg:ubar} to the 
space, say, 
\begin{equation}\label{reg:z2}
  C^0([0,T];H^{7/4}(\Omega)),
\end{equation}
which is continuously embedded into $C^0([0,T]\times \barO)$. Hence, one can repeat 
the argument in~{\bf (2)} by taking the space \eqref{reg:z2} for $z$. No modification is
required and one can see that the mapping $z \mapsto u$ is continuous from the space \eqref{reg:z2} 
to the space in \eqref{reg:u:fp}. Note that the space 
in \eqref{reg:u:fp} is continuously (though not compactly) embedded into 
the space in \eqref{reg:ubar}.
Hence, $\baru \mapsto u$ is completely continuous
because it is the composition of a compact map and a continuous one. 
Thus, to apply Schauder's theorem it just remains to choose 
a proper ball $B$ of the space in \eqref{reg:ubar} and prove that there exists 
a small time $T_1\le T$ such that the image of $B$ is contained
in $B$. This fact can be verified by a number of simple checkings. In particular,
we may use the fact that 
\begin{equation}\label{tpiccolo}
  \| v \|_{L^4(0,T_1;W^{2,3}(\Omega))} 
  \le c \| v \|_{L^4(0,T_1;H^{3}(\Omega))} 
   \le c T_1^{1/4} \| v \|_{L^\infty(0,T_1;H^3(\Omega))}
\end{equation}
for any $v\in L^\infty(0,T_1;H^3(\Omega))$. 
As a consequence, Schauder's theorem provides existence of a solution
to \eqref{eq:u:ee}-\eqref{eq:z:ee} with the initial and boundary conditions
(including \eqref{co:51}) over the time interval $(0,T_1)$. Note that, actually,
$T_1$ may be strictly smaller than $T_0$. On the other hand, performing 
the a priori estimates by keeping $\delta>0$ fixed at a first stage, we can easily see that 
the resulting bounds are uniform over the interval $(0,T)$. Hence, by
standard extension arguments, the solution to the regularized problem
can be thought to be defined over the whole of $(0,T)$. 
\beos\label{zt}
 One can see in particular that the additional regularity on $u$ obtained in
 the framework of the regularized problem is sufficient to justify the a-priori estimates
 of the previous part. Concerning $z$ there is just a point that needs to be 
 clarified a bit. Indeed, in the above part we have used the test function $-\Delta z_t$ 
 in a parabolic equation having the following structure:
 \begin{equation}\label{eql2}
   \alpha(z_t) + z_t - \Delta z \ni \eta, 
 \end{equation}
 where one can easily check that 
 \begin{equation}\label{eql2b}
   \eta = - \psi'(z) - \frac{T\dd'(z)}2 | \nabla \baru |^2 \in L^2(0,T;V).
 \end{equation}
 On the other hand, if $\alpha$ is not regularized, up to our knowledge no
 $L^2$-regularity theory is available for equation \eqref{eql2}, i.e.,
 the single summands on the \lhs\ of \eqref{eql2} are not expected to
 lie separately in $L^2$, nor it does the test function $-\Delta z_t$,
 which is then not directly admissible. To overcome this issue, one should, 
 at the step {\bf (1)}, first consider a further regularization
 of \eqref{eql2}, namely
 \begin{equation}\label{eql2la}
   \alpha\lla(z_t) + z_t - \Delta z = \eta, 
 \end{equation}
 where $\alpha\lla$ is the Yosida-approximation of $\alpha$ of order $\lambda>0$
 (cf.~\cite{Ba,Br}), and notice that \eqref{eql2la} is well-posed in $L^2$. 
 Then, one can first test \eqref{eql2la} by $-\Delta z_t$ (which is allowed thanks to better regularity
 holding for $\lambda>0$) and then take $\lambda\searrow 0$ before proceeding
 with the fixed point argument. Indeed, the obtained
 a priori bound is preserved in the limit $\lambda \searrow 0$ by semicontinuity. The details,
 based on standard convex analysis tools, are left to the 
 reader (see also~\cite[Lemma 3.10 and Proof of Theorem 3.1]{Arai} for a similar procedure).
\eddos
\beos\label{ell:par}
 It is worth observing that our choice of performing an {\sl elliptic} regularization of
 \eqref{eq:u:dd} is also motivated by the fact that a parabolic regularization
 (obtained for instance by plugging a term $\epsilon u_t$ or $-\epsilon \Delta u_t$
 in place of our $\epsilon \Delta^2 u$) would not be fully compatible with the 
 estimates of the previous section. In particular, we need to estimate (cf.~\eqref{co:20})
 the $L^2$-norm of $\nabla u$ at any {\sl fixed}\/ time $t$, and that argument seems not
 to work due to the presence of an additional term depending on $u_t$.
\eddos

\paragraph{Acknowledgments.}
G.~Akagi is supported by the Alexander von Humboldt Foundation and by the Carl Friedrich von Siemens Foundation and 
by JSPS KAKENHI Grant Number JP16H03946, JP18K18715, JP20H01812 and JP17H01095. He is also deeply grateful to the 
Helmholtz Zentrum M\"unchen and the Technishce Universit\"at M\"unchen for their kind hospitality and support during his stay in Munich. 
G.~Schimperna has been partially supported by GNAMPA (Gruppo Nazionale per l'Analisi Matematica,
la Probabilit\`a e le loro Applicazioni) of INdAM (Istituto Nazionale di Alta Matematica).

\paragraph{Conflict of Interest.}
The authors declare that they have no conflict of interest.



\vspace{15mm}

\noindent%
{\bf First author's address:}\\[1mm]
Goro Akagi\\
Mathematical Institute  and Graduate School of Science, Tohoku University,\\
6-3 Aoba, Aramaki, Aoba-ku, Sendai 980-8578 Japan\\
E-mail:~~{\tt goro.akagi@tohoku.ac.jp}

\vspace{4mm}

\noindent%
{\bf Second author's address:}\\[1mm]
Giulio Schimperna\\
Dipartimento di Matematica, Universit\`a degli Studi di Pavia\\
Via Ferrata, 5,~~I-27100 Pavia,~~Italy\\
E-mail:~~{\tt giusch04@unipv.it}

\end{document}